# Positivity of Cumulative Sums for Multi-Index Function Components Explains the Lower Bound Formula in the Levin-Robbins-Leu Family of Sequential Subset Selection Procedures

**Bruce Levin and Cheng-Shiun Leu**
Department of Biostatistics, Columbia University,
New York, New York, USA

**Abstract:** We exhibit some strong positivity properties of a certain function which implies a key inequality that in turn implies the lower bound formula for the probability of correct selection in the Levin-Robbins-Leu family of sequential subset selection procedures for binary outcomes. These properties provide a more direct and comprehensive demonstration of the key inequality than was discussed in Levin and Leu (2013).

## 1. INTRODUCTION

Levin and Leu (2013) discussed an intriguing algebraic inequality that underlies a useful lower bound formula for the probability of correct selection in the Levin-Robbins-Leu (LRL) family of sequential subset selection procedures for binomial populations. The first part of that paper presented the motivating context for the key inequality, namely, sequential subset selection procedures for selecting the best subsets of given size $b$ from $c$ competing binomial populations. LRL procedures possess useful adaptive features such as sequential elimination of apparently inferior populations and/or sequential recruitment of apparently superior populations. That the LRL family admits of a simple lower bound formula for the probability of correct select means it can easily be used to design a wide

________________________________________

Address correspondence to B. Levin, Department of Biostatistics, Mailman School of Public Health, Columbia University, 722 West 168[th] Street, New York, NY 10032, USA; E-mail: <Bruce.Levin@Columbia.edu>.

variety of selection experiments, and we showed how the key inequality implies the lower bound formula. Interested readers are encouraged to review this part of the paper for a detailed discussion of the motivating context. The second part of Levin and Leu (2013) demonstrated the validity of the key inequality assuming certain positivity properties for certain systematic components of the inequality. In the discussion section the authors pointed out that a more straightforward demonstration of the key inequality would have relied on the general positivity of certain *cumulative sums of component functions* (or *cusums* for short, not to be confused with the cusum procedures used in continuous monitoring schemes for industrial processes and other purposes). Levin and Leu (2013) established some positivity results for a limited set of cusums which were sufficient for their purpose, but they did not establish the indicated positivity results in full generality. They concluded that the positive-cusum approach was a promising strategy for a completely general proof of the key inequality, but stopped there.

The purpose of the present article is (i) to identify a wider system of component functions which might possess the positive-cusum property; and then (ii) to demonstrate how these new, stronger positivity properties follow generally from one very special case, which we state as a new lemma. We conjecture that the lemma itself can be verified for any given number of competing populations $c$ by the same computer-assisted symbolic manipulation algorithm which was discussed in Levin and Leu (2013). While this conjecture is still an open problem for *all* values of $c$, we have verified it for up to $c=6$ populations. We also offer proofs of two lemmas required for the previous approach and these new proofs offer insight into the general case. The new approach at once explains why certain single- and double-cusums had to appear in our earlier work and why the key inequality—and therefore the lower bound formula—holds.

## 2. DEFINITIONS AND NOTATION

We begin by reviewing some definitions and notation from Levin and Leu (2013). Given integers $1 \leq b < c$, an integer $k$ with $1 \leq k \leq c-b$, and constant $r > 0$, let $a=(a_1,\ldots,a_c)$ be an *ordered configuration vector*, i.e., a $c$-vector with $c-k$ positive components and $k$ zero components satisfying $a_1 \geq \cdots \geq a_b = r \geq a_{b+1} \geq \cdots \geq a_{c-k} > a_{c-k+1} = \cdots = a_c = 0$. Though the primary application to the abovementioned sequential subset selection procedures uses only integer components of $a$, here we allow the components of $a$ to be arbitrary non-negative real numbers, including the fixed component $a_b = r$. Let $a^* = (r,\ldots,r, 0,\ldots,0)$ with the first $b$ components equal to $r$. In the language of majorization (see, e.g., Marshall, et al., 2011), $a$ weakly submajorizes $a^*$ from below, written $a \succ_w a^*$. Let $B=\{1,\ldots,b\}$, $C=\{1,\ldots,c\}$, and let $\wp_C$ be the group of permutations over $C$. For any subset $K \subseteq C \setminus B$ of size $|K|=k$, let $\wp_{C \setminus K}$ be the subgroup of permutations that keep elements of $K$ fixed. Finally, let $\sigma_K$ be the permutation that sends $K$ into $Z = \{c-k+1,\ldots,c\}$ by transposing the $q^{th}$ element of $K \setminus Z$ with the $q^{th}$ element of $Z \setminus K$ for $q=1,\ldots,k-|K \cap Z|$, with $\sigma_K$ equal to the identity permutation if $K=Z$. We shall be using left cosets of the form $\sigma_K \wp_{C \backslash K}$. Permutations in such cosets are understood to apply an element of $\wp_{C \backslash K}$ first (which leaves elements in $K$ fixed), followed by $\sigma_K$.

Levin and Leu (2013) define the function $F(a)$ as follows.

$$F(a) = \sum_{\substack{K \subseteq C \backslash B \\ |K|=k}} \left\{ \frac{\sum_{\pi \in \sigma_K \wp_{C \backslash K}} w_\pi^a}{\sum_{\pi \in \sigma_K \wp_{C \backslash K}} w_\pi^{a^*}} - \frac{\sum_{\pi \in \wp_C} w_\pi^a}{\sum_{\pi \in \wp_C} w_\pi^{a^*}} \right\} = \left\{ \sum_{\substack{K \subseteq C \backslash B \\ |K|=k}} \frac{\sum_{\pi \in \sigma_K \wp_{C \backslash K}} w_\pi^a}{\sum_{\pi \in \sigma_K \wp_{C \backslash K}} w_\pi^{a^*}} \right\} - \frac{\sum_{\pi \in \wp_C} w_\pi^a}{\sum_{\pi \in \wp_C} w_\pi^{a^*}} \binom{c-b}{k}, \tag{2.1}$$

where $w$ is any given ordered vector of numbers satisfying $w_1 \geq \cdots \geq w_c \geq 1$ and where for any permutation $\pi \in \wp_C$, the symbol $w_\pi^a$ is defined as $w_\pi^a = w_{\pi(1)}^{a_1} \cdots w_{\pi(c)}^{a_c}$. $F(a^*)=0$, of course, but Theorem 4.1 of Levin and Leu (2013) asserts that $F(a) > 0$ for any $a \succ_w a^*$ with $a \neq a^*$. This is the key inequality they sought to establish. We call an expression of the form $\sum_{\pi \in \Im} w_\pi^a \Big/ \sum_{\pi \in \Im} w_\pi^{a^*}$, where the summations are over some subset of permutations $\Im \subseteq \wp_C$, a *standardized Muirhead ratio*, by extension of terminology used in Muirhead's inequality; see Muirhead (1903) and Marshall, et al. (2011). The assertion $F(a) > 0$ means that the average of *standardized Muirhead ratios* over permutation subsets $\Im$ of the form $\sigma_K \wp_{C \backslash K} \subset \wp_C$ as $K$ varies over subsets of size $k$ contained in $C \backslash B$ is bounded from below by the standardized Muirhead ratio over the full permutation group $\wp_C$.

Below we study various component sums of $F(a)$. For any $i$ in the set $I = \{1,\ldots,c-k\}$ and any $j=1,\ldots,c$, define the *single-index F-components* to be

$$F_{(j)}^{(i)}(a) = \left\{ \sum_{\substack{K \subseteq C \backslash B \\ |K|=k}} \frac{\sum_{\substack{\pi \in \sigma_K \wp_{C \backslash K} \\ \pi(i)=j}} w_\pi^a}{\sum_{\pi \in \sigma_K \wp_{C \backslash K}} w_\pi^{a^*}} \right\} - \frac{\sum_{\substack{\pi \in \wp_C \\ \pi(i)=j}} w_\pi^a}{\sum_{\pi \in \wp_C} w_\pi^{a^*}} \binom{c-b}{k} \tag{2.2}$$

and for $h=1,\ldots,c$, define the *single-index cumulative sum functions* (or *single-cusums*),

$$S_h^{(i)}(a) = \sum_{j=1}^h F_{(j)}^{(i)}(a). \tag{2.3}$$

$F_{(j)}^{(i)}(a)$ separates out all the components of $F(a)$ where $a_i$ appears as the exponent of $w_j$. Note that $F(a) = F_{(1)}^{(i)}(a) + \cdots + F_{(c)}^{(i)}(a) = S_{(c)}^{(i)}(a)$ for *any* $i \in I$.

Next, for any two superscripts $i < i' \in I$ and any two subscripts $j \neq j' \in C$ (regardless of order), define the *double-index F-components* to be

$$F_{(j,j')}^{(i,i')}(a) = \left\{ \sum_{\substack{K \subseteq C \setminus B \\ |K|=k}} \frac{\sum_{\substack{\pi \in \sigma_K \wp_{C \setminus K} \\ \pi(i)=j, \pi(i')=j'}} w_\pi^a}{\sum_{\pi \in \sigma_K \wp_{C \setminus K}} w_\pi^{a^*}} - \frac{\sum_{\substack{\pi \in \wp_C \\ \pi(i)=j, \pi(i')=j'}} w_\pi^a}{\sum_{\pi \in \wp_C} w_\pi^{a^*}} \right\} \binom{c-b}{k}. \tag{2.4}$$

Extend the notation to an arbitrary pair of subscripts $(j, j') \in C$ by setting $F_{(j,j)}^{(i,i')}(a) = 0$ whenever $j = j'$. $F_{(j,j)}^{(i,i')}(a)$ separates out all the components of $F(a)$ where $a_i$ appears as the exponent of $w_j$ and $a_{i'}$ appears as the exponent of $w_{j'}$. Note that for any $i < i' \in I$, $F_{(j)}^{(i)}(a) = \sum_{j'=1}^{c} F_{(j,j')}^{(i,i')}(a)$ and $F_{(j')}^{(i')}(a) = \sum_{j=1}^{c} F_{(j,j')}^{(i,i')}(a)$. We imagine $F_{(j,j')}^{(i,i')}(a)$ arranged in the following table (with zeros on the main diagonal).

$$\begin{array}{ccccc|c}
F_{(1,1)}^{(i,i')}(a) = 0 & \cdots & F_{(1,j')}^{(i,i')}(a) & \cdots & F_{(1,c)}^{(i,i')}(a) & F_{(1)}^{(i)}(a) \\
\vdots & \ddots & \vdots & \ddots & \vdots & \vdots \\
F_{(j,1)}^{(i,i')}(a) & \cdots & F_{(j,j')}^{(i,i')}(a) & \cdots & F_{(j,c)}^{(i,i')}(a) & F_{(j)}^{(i)}(a) \\
\vdots & \ddots & \vdots & \ddots & \vdots & \vdots \\
F_{(c,1)}^{(i,i')}(a) & \cdots & F_{(c,j')}^{(i,i')}(a) & \cdots & F_{(c,c)}^{(i,i')}(a) = 0 & F_{(c)}^{(i)}(a) \\
\hline
F_{(1)}^{(i')}(a) & \cdots & F_{(j')}^{(i')}(a) & \cdots & F_{(c)}^{(i')}(a) & F(a)
\end{array} \tag{2.5}$$

Next, for any $h, h' \in C$, define the *double cusums* $S_{(h,h')}^{(i,i')}(a) = \sum_{j=1}^{h} \sum_{j'=1}^{h'} F_{(j,j')}^{(i,i')}(a)$. The double cusums are the sums over any subrectangle of terms in the table starting at the upper left-hand corner and ending at the term $F_{(j,j)}^{(i,i')}(a)$ with $j = h$ and $j' = h'$.

In their demonstration that $F(a) > 0$, Levin and Leu (2013) established the positivity of the single cusums $S_{(h)}^{(1)}(a)$ for arbitrary $a$ and, to obtain that result, they required the positivity of the double cusums $S_{(1,h')}^{(b-1,i')}(a)$ in the limit as $a$ approaches $a^*$ through weakly submajoring configurations (see their Lemma 4.3). As a convenience, we write these limiting values as $S_{(1,h')}^{(b-1,i')}(a^+)$, where $a^+ = (r, \ldots, r, 0^+, \ldots, 0^+, 0, \ldots, 0)$ with $b$ $r$'s, $k$ zeros, and $c-b-k$ "$0^+$" entries, and where "$0^+, \ldots, 0^+$" is to be thought of as a sequence of ordered positive numbers approaching 0. The purpose of this notation is to keep the number of zero's in the sequence of configurations fixed unambiguously at the given $k$ (and not equal to the number of zeros in $a^*$, which is $c-b$) without having to include $k$ in the notation for $F$. Because the functions under consideration are continuous, the limit has a unique value irrespective of the ordered positive sequence chosen for "$0^+, \ldots, 0^+$". For example, $S_{(c)}^{(1)}(a^+) = F(a^+) = F(a^*) = 0$.

As we shall see, the positivity of cumulative sums of components of $F$ holds even more generally, namely, for multi-dimensionally indexed components, a property which here we call the *positive cusum property* (PCP) for the system of components. We define this next and state a lemma regarding the PCP for these components when evaluated at $a^+$ and then, assuming the truth of that lemma, we prove that the entire system of multidimensionally indexed components has the PCP for *every* configuration $a \succ_w a^*$ with $a \neq a^*$. This provides a more comprehensive and perhaps more elegant demonstration of the key inequality discussed in Levin and Leu (2013).

Let $p$ be an integer with $1 \leq p \leq c$, where $p$ denotes the dimension of the indexing. Let $(i) = (i_1,...,i_p)$ denote any $p$-tuple superscript with $1 \leq i_1 < \cdots < i_p \leq c-k$ and let $(j) = (j_1,...,j_p)$ denote any $p$-component vector subscript with $1 \leq j_\alpha \leq c$ for $\alpha = 1,...,p$ (irrespective of order or duplication of components). Define the *multi-index F-components* to be

$$F_{(j)}^{(i)}(a) = \left\{ \sum_{\substack{K \subseteq C \setminus B \\ |K|=k}} \frac{\sum_{\substack{\pi \in \sigma_K \wp_{C \setminus K} \\ \pi(i)=(j)}} w_\pi^a}{\sum_{\pi \in \sigma_K \wp_{C \setminus K}} w_\pi^{a^*}} \right\} - \frac{\sum_{\substack{\pi \in \wp_C \\ \pi(i)=(j)}} w_\pi^a}{\sum_{\pi \in \wp_C} w_\pi^{a^*}} \binom{c-b}{k}, \tag{2.6}$$

which appears almost identical to (2.2) but where the notation $\pi(i) = (j)$ refers to the restriction of permutations in the numerators to those for which $\pi(i_\alpha) = j_\alpha$ for each $\alpha = 1,...,p$. As above, if any pair of subscripts $j_\alpha = j_{\alpha'}$ for some $\alpha \neq \alpha'$, set $F_{(j)}^{(i)}(a) = 0$. Summation of $F_{(j)}^{(i)}(a)$ over any subset of subscript elements in $(j)$ from 1 to $c$ reduces to a lower-dimension multi-index $F$-component with corresponding superscript elements removed from $(i)$. Finally, let the $p$-vector $(h) = (h_1,...,h_p)$ with $1 \leq h_\alpha \leq c$ for $\alpha = 1,...,p$ (irrespective of order or duplication of components) be an arbitrary multi-index subscript over which to cumulatively sum the multi-components of $F$ and define the *multi-cusum*

$$S_{(h)}^{(i)}(a) = \sum_{j_1=1}^{h_1} \cdots \sum_{j_p=1}^{h_p} F_{(j)}^{(i)}(a). \tag{2.7}$$

We know in advance that certain multi-cusums must be zero, namely, those in which *every* term in the sum (2.7) involves a multi-index $F$ component with repeated elements in its subscript. For example, suppose $b=2$, $c=4$, $a=(a_1, r, a_3, 0)$ with $k=1$, and $p=3$, in which case there is only one multi-index superscript to consider, $(i)=(1,2,3)$. For $(h)=(1, 2, 2)$, the multi-index $F$ components in (2.7) have subscripts $(1,1,1)$, $(1,1,2)$, $(1,2,1)$, and $(1,2,2)$, every one of which has a subscript with repeated elements. Thus, we know in advance that such a multi-cusum must be zero. Therefore, we define a *contributing* cusum subscript $(h)$ as one in which *at least one* multi-index $F$-component in (2.7) has a multi-index subscript $(j)$ *with all elements distinct*, $j_\alpha \neq j_{\alpha'}$ for all pairs $\alpha, \alpha'$ with $1 \leq \alpha \neq \alpha' \leq p$. It is easy to see that a cusum subscript $(h) = (h_1,...,h_p)$ is contributing if and only if after rearranging the elements in ascending order as, say, $h'_1 \leq \cdots \leq h'_p$ one has $h'_\alpha \geq \alpha$ for each $\alpha = 1,\ldots,p$. We define a *contributing multi-cusum* as one with a contributing cusum subscript $(h)$.

The following lemma asserts the PCP for all contributing multi-cusums evaluated at the special configuration $a = a^+$.

**Lemma 2.1.** *For any $p=1,\ldots,c-k$, any p-tuple $(i) = (i_1,\ldots,i_p)$ with $1 \leq i_1 < \cdots < i_p \leq c-k$, and any contributing cusum subscript $h = (h_1,\ldots,h_p)$ with $1 \leq h_\alpha \leq c-1$ for $\alpha = 1,\ldots p$, we have $S^{(i)}_{(h)}(a^+) > 0$, i.e., all contributing multi-cusums are positive when evaluated at $a^+ = (r,\ldots,r,0+,\ldots,0+,0,\ldots,0)$.*

Remarkably, exactly the same algorithm for computer-assisted symbolic manipulation that was stated in Appendix B of Levin and Leu (2013a) suffices to prove Lemma 2.1 in all cases checked to date, comprising thousands of scenarios with different $b$, $c$, $k$, $p$, $(i)$, and $(h)$. The complexity of the calculations grows exponentially, though, so we stopped with $c=6$, but we conjecture that the algorithm would continue to suffice for any larger $c$. We further analyze Lemma 2.1 in Section 4.

For convenience in what follows, instead of limiting the cusum subscripts to values no greater than $c-1$, we shall proceed equivalently by fixing $p=c-k$ and allowing any cusum subscript element to take the value $c$. As mentioned above (2.7), summation from 1 to $c$ over any subscript element or subset of elements of $F^{(i)}_{(j)}(a)$ reduces $S^{(i)}_{(h)}(a)$ to a lower-order multi-cusum with corresponding superscript elements removed from $(i)$ and with all remaining subscripts constrained to the range $1,\ldots,c-1$, as in the Lemma. Doing so will allow us to suppress the superscripts from the notation, which henceforth are assumed to be fixed at the value $(i)=(1,\ldots,c-k)$, and all cusum subscripts will be $(c-k)$-vectors in the integer lattice $C^{c-k}$. Also, for any $i=1,\ldots,c-k$, it will be convenient to rearrange the order of summation in the multi-cusum and abbreviate iterated summations as follows. We write

$$S_{(h)}(a) = \sum_{j_1=1}^{h_1}\cdots\sum_{j_p=1}^{h_p} F_{(j)}(a) = \sum_{j_i=1}^{h_i}\sum_{j_1=1}^{h_1}\cdots\sum_{j_{i-1}=1}^{h_{i-1}}\sum_{j_{i+1}=1}^{h_{i+1}}\cdots\sum_{j_p=1}^{h_p} F_{(j)}(a) = \sum_{j_i=1}^{h_i}\sum_{j'_i=1}^{h'_i} F_{(j_i,j'_i)}(a) = S_{(h_i,h'_i)}(a),$$

where $j'_i = (j_1,\ldots,j_{i-1},j_{i+1},\ldots,j_p)$, $h'_i = (h_1,\ldots,h_{i-1},h_{i+1},\ldots,h_p)$, the inner iterated summations are abbreviated as $\sum_{j'_i=1}^{h'_i} F_{(j_i,j'_i)}(a) = \sum_{j_1=1}^{h_1}\cdots\sum_{j_{i-1}=1}^{h_{i-1}}\sum_{j_{i+1}=1}^{h_{i+1}}\cdots\sum_{j_p=1}^{h_p} F_{(j_i,j'_i)}(a)$, and where the split subscripts are intended to denote $(j)$ or $(h)$ in their proper order, i.e., $(j_i,j'_i) = (j) = (j_1,\ldots,j_i,\ldots,j_p)$ or $(h_i,h'_i) = (h) = (h_1,\ldots,h_i,\ldots,h_p)$, respectively. With these conventions, we may now prove the following result.

## 3. MAIN RESULT

**Theorem 3.1.** *All contributing multi-cusums are strictly positive. More precisely, if (h) is a contributing cusum subscript with $1 \leq h_\alpha \leq c$ for $\alpha = 1,...,c-k$, then $S_{(h)}(a) = S_{(h)}^{(1,...,c-k)}(a) > 0$ for every ordered configuration $a \succ_w a^*$ with $a \neq a^*$.*

*Proof.* For any contributing subscript $(h)$, define the function $G(x) = G_{(h)}(x) = S_{(h)}(r,...,r,x,...,x,0,...0)$ with the first $b$ elements equal to $r$, the last $k$ elements equal to zero, and the middle $c-b-k$ elements equal to $x$ for $0 < x \leq r$. It will be convenient to abbreviate such configurations as "$\cdot x \cdot$", writing, e.g., $G(x) = S_{(h)}(\cdot x \cdot)$. Extend the definition of $G$ to $x=0$ by continuity as $G(0^+)$. Observe that $G(0^+) = S_{(h)}(a^+) > 0$ by Lemma 2.1. We shall be interested in the derivative of $G$ with respect to $x$, which involves partial derivatives of multi-index components of $F$ with respect to components of $a$. To that end, we note from definition (2.6) that we have

$$(\partial/\partial a_i)F_{(j)}(a) = (\partial/\partial a_i)F_{(j_1,...,j_{c-k})}(a) = (\partial/\partial a_i)w_{j_i}^{a_i} \times \{\text{terms not involving } a_i\}$$

$$= (\log w_{j_i})w_{j_i}^{a_i} \times \{\text{the same terms not involving } a_i\} = F_{(j)}(a)\log w_{j_i}.$$

Then

$$(d/dx)G(0^+) = \sum_{i=b+1}^{c-k}(\partial/\partial a_i)S_{(h)}(a^+)$$

and we have the following key chain of inequalities for the partial derivatives.

$$(\partial/\partial a_i)S_{(h)}(a^+) = \sum_{j_i=1}^{h_i}\sum_{j'_i=1}^{h'_i}(\partial/\partial a_i)F_{(j_i,j'_i)}(a^+)$$

$$= \sum_{j_i=1}^{h_i}\sum_{j'_i=1}^{h'_i}F_{(j_i,j'_i)}(a^+)\log w_{j_i}$$

$$= \sum_{j'_i=1}^{h'_i}F_{(1,j'_i)}(a^+)\log w_1 + \sum_{j'_i=1}^{h'_i}F_{(2,j'_i)}(a^+)\log w_2 + \cdots + \sum_{j'_i=1}^{h'_i}F_{(h_i,j'_i)}(a^+)\log w_{h_i}$$

$$= S_{(1,h'_i)}(a^+)\log w_1 + \sum_{j'_i=1}^{h'_i}F_{(2,j'_i)}(a^+)\log w_2 + \cdots + \sum_{j'_i=1}^{h'_i}F_{(h_i,j'_i)}(a^+)\log w_{h_i}$$

$$\geq S_{(1,h'_i)}(a^+)\log w_2 + \sum_{j'_i=1}^{h'_i}F_{(2,j'_i)}(a^+)\log w_2 + \cdots + \sum_{j'_i=1}^{h'_i}F_{(h_i,j'_i)}(a^+)\log w_{h_i}$$

$$= \left\{S_{(1,h'_i)}(a^+) + \sum_{j'_i=1}^{h'_i}F_{(2,j'_i)}(a^+)\right\}\log w_2 + \sum_{j'_i=1}^{h'_i}F_{(3,j'_i)}(a^+)\log w_3 + \cdots + \sum_{j'_i=1}^{h'_i}F_{(h_i,j'_i)}(a^+)\log w_{h_i} \quad (3.1)$$

$$= S_{(2,h'_i)}(a^+)\log w_2 + \sum_{j'_i=1}^{h'_i}F_{(3,j'_i)}(a^+)\log w_3 + \cdots + \sum_{j'_i=1}^{h'_i}F_{(h_i,j'_i)}(a^+)\log w_{h_i}$$

$$\geq S_{(2,h'_i)}(a^+)\log w_3 + \sum_{j'_i=1}^{h'_i}F_{(3,j'_i)}(a^+)\log w_3 + \cdots + \sum_{j'_i=1}^{h'_i}F_{(h_i,j'_i)}(a^+)\log w_{h_i}$$

$$= \left\{S_{(2,h'_i)}(a^+) + \sum_{j'_i=1}^{h'_i}F_{(3,j'_i)}(a^+)\right\}\log w_3 + \cdots + \sum_{j'_i=1}^{h'_i}F_{(h_i,j'_i)}(a^+)\log w_{h_i}$$

$$= S_{(3,h'_i)}(a^+)\log w_3 + \cdots + \sum_{j'_i=1}^{h'_i}F_{(h_i,j'_i)}(a^+)\log w_{h_i}$$

$$\vdots$$

$$\geq S_{(d_i,h'_i)}(a^+)\log w_d$$

where $d_i = \min(\lambda, h_i)$ and $\lambda$ is the index of the "last" positive log odds, i.e., $w_{\lambda-1} > 1$ but $w_{\lambda+1} = \cdots = w_c = 1$. To avoid trivialities we assume not all odds parameters are equal, so $1 \leq \lambda \leq c$. The inequalities hold due to the positivity of the contributing cusums stated in the Lemma. Any inequality in the chain will be strict if the cusum on the left of the inequality is contributing (in which case it is strictly positive), say $S_{(\alpha_i,h'_i)} > 0$, *and* the two odds parameters surrounding the inequality are unequal ($w_\alpha > w_{\alpha+1}$). Because not all of the sub-cusums $S_{(1,h'_i)},\ldots,S_{(d_i,h'_i)}$ are necessarily contributing and some of the odds parameters may be equal, it is possible for a partial derivative of $S_{(h)}(a)$ to equal zero at $a=a^+$. If this situation pertains to each of the partial derivatives, it is even possible for the derivative of $G(x)$ to equal zero at $x=0^+$. In such cases, however, the derivative will be zero *for all x*, because a non-contributing cusum index is only a feature of the indexing, not the argument $x$, and similarly $\lambda$ depends only on the odds parameters. Thus $G(x)$ would be constant and equal to $S_{(h)}(a^+) > 0$ by Lemma 2.1. For example, suppose $b=2$, $c=4$, $k=1$, and consider the cusum subscript $(h)=(2,1,3)$ which is contributing as it allows one multi-index component subscript, namely, $(j)=(2,1,3)$ with distinct elements. However, if $w_3 = w_4 = 1$, the partial derivative of $S_{(h)}(a)$ with respect to $a_3$ at $a^+$ actually equals $(\partial/\partial a_3) S_{(h)}(a^+) = (\partial/\partial a_3) F_{(2,1,3)}(a^+) = F_{(2,1,3)}(a^+) \log w_3 = 0$. This is consistent with the end of the inequality chain (3.1) because $d_3 = \min(\lambda, 3) = 2$ and so $S_{(2,1,2)}(a^+) \log w_2 = 0$ because (2,1,2) is not a contributing cusum index. Thus $G(x) = S_{(2,1,3)}(r, r, x, 0)$ is constant and equal to $S_{(h)}(a^+) > 0$.

In all other cases there will be at least one partial derivative for which at least one of the inequalities in (3.1) is strict because unequal odds consecutively multiply a positive contributing multi-cusum, whereupon (3.1) becomes $(\partial/\partial a_i) S_{(h)}(a^+) > S_{(d_i, h_i')}(a^+) \log w_{d_i}$ with a strict inequality. Thus, even if the cusum on the right-hand side were itself non-contributing, the derivative of $G$ at $0^+$ would still be positive. As a result, even if $G$ itself were zero at $0^+$, e.g., if $(h)=(c,\ldots,c)$, $G(x)$ would nevertheless increase in a neighborhood to the right of 0.

Therefore, in all cases there exists a $\delta > 0$ such that $G(x) > 0$ for all $x$ in the open interval $0 < x < \delta$. Furthermore, we can choose $\delta$ such that $S_{(\alpha_i, h_i')}(x) > 0$ in the interval for *every* contributing multi-cusum with subscript $1 \leq \alpha_i \leq h_i$ for each $i = b+1, \ldots, c-k$, because there are only finitely many of these. We claim that in fact $G(x) > 0$ *for all* $0 < x \leq r$. We prove this next by contradiction.

Suppose to the contrary that there exists an $x'$ with $0 < x' \leq r$ such that $G(x') \leq 0$. Consider the derivative of $G$ at $x'$ and the chain of inequalities in (3.1) applied at the argument $\cdot x' \cdot$ rather than $a^+$ in the cusums. At least one contributing cusum must be non-positive somewhere to the left of $x'$, i.e., there must exist an $x''$ with $0 < x'' < x'$ such that $S_{(\alpha_i, h_i')}(\cdot x'' \cdot) \leq 0$ for some $i = b+1, \ldots, c-k$ and some $1 \leq \alpha_i \leq d_i \leq h_i$. This is so because if each contributing sub-cusum were positive everywhere in the open interval $(0, x')$, then the chain of inequalities in (3.1) evaluated at any $x$ in the interval would imply that $G(x)$ is positive and increasing (or at worst constant) throughout the interval and therefore at $x'$ by continuity, contradicting the assumption that $G(x') \leq 0$. Then for each contributing sub-cusum $S_{(\alpha_i, h_i')}$ which is somewhere non-positive, the corresponding set $M(\alpha_i, h_i') = \{x : \delta \leq x \leq x' \text{ and } S_{(\alpha_i, h_i')}(\cdot x \cdot) \leq 0\}$ is non-empty and each one is obviously bounded from below by $\delta$. Let $x^*$ be the minimum of their greatest lower bounds, i.e., $x^* = \min_{i, \alpha_i} \{\inf x : x \in M(\alpha_i, h_i')\}$, with $\delta \leq x^* \leq x'$. Suppose the minimum occurs for the set $M(\beta_{i*}, h_{i*}')$. The corresponding contributing cusum must equal zero at $x^*$, $S_{(\beta_{i*}, h_{i*}')}(x^*) = 0$. This is because for values to the left of $x^*$, the cusum must be positive, else $x^*$ would not be a lower bound for the corresponding set $M(\beta_{i*}, h_{i*}')$; and to the right of $x^*$, there are values of $x$ arbitrarily close to $x^*$ at which the cusum is negative, else $x^*$ would not be the *greatest* lower bound for $M(\beta_{i^*}, h'_{i^*})$. Then by continuity, we have $S_{(\beta_{i^*}, h'_{i^*})}(x^*) = 0$. (Even if $x^* = x' = r$ with no value of $x$ to the right of $x^*$, $S_{(\beta_{i^*}, h'_{i^*})}(x^*)$ still would equal 0 by continuity.)

We are now ready to derive a contradiction. There are two possibilities to consider. One is that the identified multi-cusum $S_{(\beta_{i^*}, h'_{i^*})}$ is actually the original cusum $S_{(h)}$ which defines $G(x) = S_{(h)}(\cdot x \cdot)$; this would be the case if $d_{i^*} = \min(\lambda, h_{i^*}) = h_{i^*}$ and $\inf M(h_{i^*}, h'_{i^*})$ is the minimum infimum among all the non-empty sets $M(\alpha_i, h'_i)$. The second possibility is that the identified contributing cusum $S_{(\beta_{i^*}, h'_{i^*})}$ is a strict sub-cusum, i.e., $\beta_{i^*} < h_{i^*}$.

In the first case, by our choice of $\delta$, there are positive values of $x$ near zero, say at $x_0$ with $0 < x_0 < \delta$, such that $G(x_0) > 0$ and thus the slope of the straight line connecting $(x_0, G(x_0))$ and $(x^*, G(x^*)) = (x^*, 0)$ is negative. By the mean value theorem for derivatives, there exists an $x'''$ with $x_0 < x''' < x^*$ such that $(d/dx)G(x''') < 0$. But *all* contributing cusums in the chain of inequalties (3.1) are positive when evaluated at any point to the left of $x^*$, because $x^*$ was chosen as the smallest of the infimums for non-empty sets $M(\alpha_i, h'_i)$ where there were possible violations of the chain of inequalities; and other contributing cusums are positive everywhere in $(\delta, x')$. Hence by (3.1) evaluated at $x'''$, $(d/dx)G(x''') \geq 0$, contradiction. Therefore, in this first case, $G(x)$ must be positive for all $x$ in the interval $(0, r]$.

In the second case, we simply repeat the preceding argument for the contributing sub-cusum $S_{(\beta_{i^*}, h'_{i^*})}$, redefining $G(x)$ to be $G(x) = S_{(\beta_{i^*}, h'_{i^*})}(\cdot x \cdot)$ with $G(x^*) = 0$. By an obvious inductive argument, eventually we must arrive at a contributing cusum violating the chain of inequalities that would fall into the first case, which would lead to a contradiction as above. For example, consider $G(x) = S_{(c,\ldots,c,h)}(\cdot x \cdot)$, which is the single-index cusum $S_{(h)}^{(c-k)}(a)$ evaluated at $a = (r,\ldots,r,x,\ldots,x,0,\ldots,0)$. If all cusums are positive in (3.1) evaluated at a given $x$, then $(\partial/\partial a_{c-k})S_{(h,c)}(\cdot x \cdot) \geq S_{(d,c)}(\cdot x \cdot) \log w_d$, where $d = \min(\lambda, h)$ while for $i = b+1, \ldots, c-k-1$, $(\partial/\partial a_i)S_{(c,h'_i)}(\cdot x \cdot) \geq S_{(d_i, h'_i)}(\cdot x \cdot) \log w_{d_i} = S_{(\lambda, h')}(\cdot x \cdot) \log w_\lambda$, where $h' = (c,\ldots,c,h)$. Each of the latter terms are equal because the cusums are identical when evaluated at the replicated elements in $\cdot x \cdot$. Therefore, in single- and double-cusum notation, $(d/dx)G(x) \geq S_{(d)}^{(c-k)}(\cdot x \cdot) \log w_d + (c-b-k-1)S_{(\lambda,d)}^{(i,c-k)}(\cdot x \cdot) \log w_\lambda$. If, on the contrary, $G(x')$ were non-positive for some $\delta < x' \leq r$, then at some $0 < x'' < x'$ either $S_{(q)}^{(c-k)}(\cdot x'' \cdot) \leq 0$ for some $q \leq d \leq h$ or $S_{(q',q)}^{(i,c-k)}(\cdot x'' \cdot) \leq 0$ for some $q' \leq \lambda$. If $q = h$ and $q' = \lambda = c$, the first case obtains. If $q < h$ or $q' < c$, then the second case obtains. In the first case, one proceeds to a contradiction as above. In the second case, we (inductively) apply the *reductio ad absurdum* argument to the sub-cusum $S_{(q)}^{(c-k)}(\cdot x'' \cdot) \leq 0$ or $S_{(q',q)}^{(i,c-k)}(\cdot x'' \cdot) \leq 0$. Eventually the first case must obtain.

Therefore, in all cases, we have shown that $G(x)$ is positive on $(0, r]$ and strictly increases throughout (or at worst is a positive constant). The proof furthermore showed that all contributing sub-cusums of the form $S_{(\alpha_i, h'_i)}(\cdot x \cdot)$ are positive for each $i = b+1, \ldots, c-k$.

To complete the proof, we simply repeat all of the above arguments for $G(x)$ redefined next as $G(x) = S_{(h)}(r, \ldots, r, x, \ldots, x, a_{c-k}, 0, \ldots, 0)$ for $x$ in the interval $[a_{c-k}, r]$, noting that $G(a_{c-k}) > 0$ by the preceding step. Continuing toward the left in this way, we have that all contributing cusums of the form $S_{(h)}(r, \ldots, r, a_{b+1}, \ldots, a_{c-k}, 0, \ldots, 0)$ are positive. Then we repeat the argument with $G(x)$ redefined as $G(x) = S_{(h)}(x, \ldots, x, r, a_{b+1}, \ldots, a_{c-k}, 0, \ldots, 0)$ with $b-1$ $x$'s for $x \geq r$, followed by the same argument with $G(x) = S_{(h)}(x, \ldots, x, a_{b-1}, r, x, \ldots, x, a_{c-k}, 0, \ldots, 0)$ for $x \geq a_{b-1}$, and so on. We conclude finally that $S_{(h)}(a) > 0$ for any contributing cusum index and any $a \succ_W a^*$ with $a \neq a^*$. Q.E.D.

## 4. FURTHER REMARKS ON LEMMA 2.1

Beyond the computational confirmations of Lemma 2.1 produced for $c \leq 6$, there are several illuminating special cases in which the lemma can be proven to hold for all $c$ and which provide heuristic reasons why it appears to hold in complete generality. To explain them, let us re-express (2.7) using *tuplet notation*, say $(d)$, which is a convenient way to denote a generic $d$-tuplet $(d) = (i_1, \ldots, i_d)$ with $1 \leq i_1 < \cdots < i_d \leq c$, and let $w_{(d)}$ denote the product $w_{(d)} = w_{i_1} \cdots w_{i_d}$. To avoid confusion with the previous notation for multi-index super- or subscripts, we will use square brackets for the multi-index superscripts and subscripts, e.g., $S^{[i]}_{[h]}(a)$. Furthermore, when specifying set operations involving the subset of multi-index elements (i.e., without regard to order), we will use the usual curly brackets of set notation, e.g., $(C \setminus B) \setminus \{j\}$. There are three facts about $F$-components and multi-cusums evaluated at $a = a^+$.

(i) If any $K \subseteq C \setminus B$ of size $|K| = k$ contains any element of subscript $[j]$, then $\sigma_K \wp_{C \setminus K}$ has no permutations that satisfy the constraints $\pi[i] = [j]$ in the left-hand numerators. Therefore, we can write the outer sum in the term in braces as a sum over subsets $K$ of size $k$ in $(C \setminus B) \setminus \{j\}$.

(ii) Odds parameters only appear in the functions evaluated at $a^+$ when associated with superscript elements $1, \ldots, b$, and when it appears, $w_j$ only appears as $w^r_j$. We may therefore take $r = 1$ without loss of generality and write $w^{a^+}_\pi = w_{\pi(1)} \cdots w_{\pi(b)} = w_{(b)}$ where the $b$-tuple $(b)$ comprises the permutation values $\pi(1), \ldots, \pi(b)$ after rearranging them in ascending order. Similarly, we will write $w^{a^*}_\pi = w_{(b)}$.

(iii) As we shall see below, when evaluating the $F$-components and multi-cusums at configuration $a^+$, the values turn out to depend on the actual superscript elements in $[i]$ only through the number of elements less than or equal to $b$. Accordingly, let $q = q[i] = \#\{\alpha = 1,...,p : i_\alpha \leq b\}$ denote that number, such that $i_q \leq b < i_{q+1}$. The upper range of $q$ is clearly $q \leq \min(b, p)$ and because there can be no more than $c - k - b$ elements of $[i]$ beyond $i_q$, we must have $p - q \leq c - k - b$ or $q \geq \min(0, b + k + p - c)$. The only odds parameters fixed by the constraint $\pi[i] = [j]$ that actually appear in the expressions have subscripts $j_1,..., j_q$. Thus $w_\pi^{a^+} = w_{\pi(1)} \cdots w_{\pi(b)} = w_{(b)} = w_{j_1} \cdots w_{j_q} \cdot w_{(b-q)}$ where $(b-q)$ is the $(b-q)$-tuple $(b - q) = (\pi(q+1),...,\pi(b)) = (j_{q+1},..., j_b)$.

We may now replace terms $w_\pi^{a^+}$ and $w_\pi^{a^*}$ with corresponding tuplet terms but in so doing we need to track the multiplicites with which different permutations $\pi$ give rise to the same tuplet. In particular,

$$F_{[j]}^{[i]}(a^+)/w_{j_1} \cdots w_{j_q} = \left\{ \sum_{\substack{K \subseteq (C \setminus B) \setminus \{j\} \\ |K| = k}} \frac{M_{left,num} \cdot \sum_{(b-q) \subseteq (C \setminus K) \setminus \{j\}} w_{(b-q)}}{M_{left,den} \cdot \sum_{(b) \subseteq C \setminus K} w_{(b)}} \right\} \cdot \frac{M_{right,num} \cdot \sum_{(b-q) \subseteq C \setminus \{j\}} w_{(b-q)}}{M_{right,den} \cdot \sum_{(b) \subseteq C} w_{(b)}} \binom{c-b}{k}$$

where the multiplicities are

$$M_{right,den} = b!(c-b)!$$

$$M_{left,den} = b!(c-b-k)!$$

$$M_{right,num} = (b-q)!(c-b-p+q)!, \text{ and}$$

$$M_{left,num} = (b-q)!(c-b-k-p+q)! \ .$$

Writing $\binom{c-b}{k} = \frac{(c-b)!}{k!(c-b-k)!}$ and rearranging the factorials, we arrive at the first of several useful re-expressions:

$$p! \left\{ \binom{b}{q}\binom{c-b-k}{p-q} \middle/ \binom{p}{q} \right\} F_{[j]}^{[i]}(a^+)/w_{j_1} \cdots w_{j_q}$$

$$= \left\{ \sum_{\substack{K \subseteq (C \setminus B) \setminus \{j\} \\ |K|=k}} \frac{\sum_{(b-q) \subseteq (C \setminus K) \setminus \{j\}} w_{(b-q)}}{\sum_{(b) \subseteq C \setminus K} w_{(b)}} \right\} \cdot \frac{\sum_{(b-q) \subseteq C \setminus \{j\}} w_{(b-q)}}{\sum_{(b) \subseteq C} w_{(b)}} \binom{c-b-p+q}{k}$$

(4.1)

and

$$p!\left\{\binom{b}{q}\binom{c-b-k}{p-q}\Big/\binom{p}{q}\right\}S_{[h]}^{[i]}(a^+)$$

$$=\sum_{j_1=1}^{h_1}\cdots\sum_{j_q=1}^{h_q}w_{j_1}\cdots w_{j_q}\cdot\sum_{j_{q+1}=1}^{h_{q+1}}\cdots\sum_{j_p=1}^{h_p}\left[\left\{\sum_{\substack{K\subseteq(C\setminus B)\setminus\{j\}\\|K|=k}}\frac{\sum_{(b-q)\subseteq(C\setminus K)\setminus\{j\}}w_{(b-q)}}{\sum_{(b)\subseteq C\setminus K}w_{(b)}}\right\}-\frac{\sum_{(b-q)\subseteq C\setminus\{j\}}w_{(b-q)}}{\sum_{(b)\subseteq C}w_{(b)}}\binom{c-b-p+q}{k}\right].$$

We will abbreviate the latter expression as

$$p!\left\{\binom{b}{q}\binom{c-b-k}{p-q}\Big/\binom{p}{q}\right\}S_{[h]}^{[i]}(a^+)$$

$$=\sum_{[j]\leq[h]}w_{(q\uparrow j)}\left\{\sum_{\substack{K\subseteq(C\setminus B)\setminus\{j\}\\|K|=k}}\frac{\sum_{(b-q)\subseteq(C\setminus K)\setminus\{j\}}w_{(b-q)}}{\sum_{(b)\subseteq C\setminus K}w_{(b)}}\right\}-\sum_{[j]\leq[h]}w_{(q\uparrow j)}\frac{\sum_{(b-q)\subseteq C\setminus\{j\}}w_{(b-q)}}{\sum_{(b)\subseteq C}w_{(b)}}\binom{c-b-p+q}{k},\quad(4.2)$$

where we have borrowed the "take" operator from the APL programming language to write the $(q)$-tuple $(j_1,\ldots,j_q)$ as $(q\uparrow j)$ and $w_{(q\uparrow j)}=w_{j_1}\cdots w_{j_q}$.

**Remark 4.1.** As noted above, the expression for the multi-cusum depends on the multi-index superscript $[i]$ only through $q=q[i]$. Therefore, when checking Lemma 2.1 using computer-assisted symbolic manipulation, one need only enumerate the various values of $c$, $b$, $k$, $p$, and $q$ and check the positivity of (4.2). If desired for concreteness, one may select a representative superscript for given $b$, $c$, $k$, $p$, and $q$ such as $[i]=(1,\ldots,q,b+1,\ldots,b+p-q)$. The relevant ranges are $1\leq b\leq c-1$, $1\leq k\leq c-b-1$ (the case $k=c-b$ being trivial), $1\leq p\leq c-k$, and $\max(0,b+k+p-c)\leq q\leq\min(b,p)$.

Even more suggestive expressions result when we replace sums with averages. To that end, we state the following additional facts.

(iv) As noted in Section 2, to be a contributing cusum subscript $[h]$ must satisfy $h'_\alpha\geq\alpha$, where the $h'_\alpha$ are the ordered cusum subscripts with $h'_1\leq\cdots\leq h'_p$ for $\alpha=1,\ldots,p$. For a given cusum subscript $[h]$, the quantity $h'_1(h'_2-1)\cdots(h'_p-p+1)$ gives the number of multi-index $F$-component subscripts $[j]$ involved in the cumulative sum $S_{[h]}^{[i]}(a^+)$. In the special case $h_1=\cdots=h_p=h$, that number would be $p!\binom{h}{p}$, so it will be convenient to extend the binomial coefficient notation by defining $\binom{[h]}{p}=h'_1(h'_2-1)\cdots(h'_p-p+1)/p!$. Then the number of $[j]\leq[h]$ is $p!\binom{[h]}{p}$.

(v) The number of (b–q)-tuples contained in $(C \setminus K) \setminus \{j\}$ is $\binom{c-k-p}{b-q}$; the number of (b–q)-tuples in $C \setminus \{j\}$ is $\binom{c-p}{b-q}$; the number of (b)-tuples in $C \setminus K$ is $\binom{c-k}{b}$; and the number of (b)-tuples in $C$ is $\binom{c}{b}$.

(vi) For a given multi-index subscript $[j] = (j_1, \ldots, j_p)$, let the subset of elements that comprise $[j]$ be denoted by $J = \{j\} = \{j_1, \ldots, j_p\}$. Then the number of subsets $K \subseteq C \setminus B$ of size $|K| = k$ that do not contain any elements of $[j]$ only depends on $J$ and is given by $k_J = \binom{c-b-p+q'(J)}{k}$, where $q'(J) = \#\{\alpha = 1, \ldots, p : j_\alpha \le b\} = q'[j]$ is the number of elements of $J$ (or $[j]$) less than or equal to $b$. Like $q = q[i]$, the range of $q'(J)$ is $\max(0, b+k+p-c) \le q'[j] \le \min(b, p)$. Unlike $q[i]$, however, the distributional properties of $q'[j] = q'(J)$ play a crucial role in cusum positivity, as explained below. For notational simplicity we shall continue to abbreviate $q[i]$ as $q$ without further comment and will write $q'$ or $q'(J)$ for clarity as the context requires.

(vii) For given multi-index superscript $[i]$ and subscript $[j]$, let $q = q[i]$ and $J = \{j\}$ as defined above. Then for given $K \subseteq C \setminus B$ with $|K| = k$ and $J \cap K = \emptyset$, define the *cross-product ratio* $R_{JK}^{(q)}(w)$ as

$$R_{JK}^{(q)}(w) = \frac{\left\{w_{(q\uparrow j)} \sum_{(b-q) \subseteq C \setminus (J \cup K)} w_{(b-q)} \Big/ \binom{c-k-p}{b-q}\right\}\left\{\sum_{(b) \subseteq C} w_{(b)} \Big/ \binom{c}{b}\right\}}{\left\{w_{(q\uparrow j)} \sum_{(b-q) \subseteq C \setminus J} w_{(b-q)} \Big/ \binom{c-p}{b-q}\right\}\left\{\sum_{(b) \subseteq C \setminus K} w_{(b)} \Big/ \binom{c-k}{b}\right\}},$$

Because $w_{(q\uparrow j)}$ appears in both numerator and denominator of $R_{JK}$, it can be written equivalently as

$$R_{JK}^{(q)}(w) = \frac{\left\{\sum_{(b-q) \subseteq C \setminus (J \cup K)} w_{(b-q)} \Big/ \binom{c-k-p}{b-q}\right\}\left\{\sum_{(b) \subseteq C} w_{(b)} \Big/ \binom{c}{b}\right\}}{\left\{\sum_{(b-q) \subseteq C \setminus J} w_{(b-q)} \Big/ \binom{c-p}{b-q}\right\}\left\{\sum_{(b) \subseteq C \setminus K} w_{(b)} \Big/ \binom{c-k}{b}\right\}}, \tag{4.3}$$

and note that $R_{JK}^{(q)}(w)$ depends on the specific values of $w_{j_1}, \ldots, w_{j_q}$ only through their number $q = q[i]$ and thus it depends on $[j]$ only through the subset $J = \{j\}$ (without regard to order). The cross-product ratios $R_{JK}$ will play a key role in the sequel.

Now the cusum in (4.2) is positive if and only if

$$\frac{\displaystyle\sum_{[j]\le[h]} w_{(q\uparrow j)} \left\{ \sum_{\substack{K\subseteq(C\setminus B)\setminus\{j\} \\ |K|=k}} \frac{\displaystyle\sum_{(b-q)\subseteq(C\setminus K)\setminus\{j\}} w_{(b-q)}}{\displaystyle\sum_{(b)\subseteq C\setminus K} w_{(b)}} \right\}}{\displaystyle\sum_{[j]\le[h]} w_{(q\uparrow j)} \frac{\displaystyle\sum_{(b-q)\subseteq C\setminus\{j\}} w_{(b-q)}}{\displaystyle\sum_{(b)\subseteq C} w_{(b)}}} > \binom{c-b-p+q}{k}$$

But the ratio on the left-hand side is equal to

$$\frac{\binom{c-k-p}{b-q}\binom{c}{b}}{\binom{c-p}{b-q}\binom{c-k}{b}} \sum_{[j]\le[h]} \left\{ \frac{w_{(q\uparrow j)} \displaystyle\sum_{(b-q)\subseteq C\setminus J} w_{(b-q)}}{\displaystyle\sum_{[j]\le[h]} w_{(q\uparrow j)} \sum_{(b-q)\subseteq C\setminus J} w_{(b-q)}} \right\} \left\{ \sum_{\substack{K\subseteq(C\setminus B)\setminus J \\ |K|=k}} R^{(q)}_{JK}(w) \right\}$$

$$= \frac{\binom{c-k-p}{b-q}\binom{c}{b}}{\binom{c-p}{b-q}\binom{c-k}{b}} \frac{\displaystyle\sum_{[j]\le[h]} w_{(q\uparrow j)} u^{(q)}_J(w) \{k_J R^{(q)}_{J\bullet}(w)\}}{\displaystyle\sum_{[j]\le[h]} w_{(q\uparrow j)} u^{(q)}_J(w)}$$

$$= \frac{\binom{c-k-p}{b-q}\binom{c}{b}}{\binom{c-p}{b-q}\binom{c-k}{b}} \frac{\displaystyle\sum_{J\subseteq C} f_J[h] \overline{w}_{(q\uparrow J)} u^{(q)}_J(w) \{k_J R^{(q)}_{J\bullet}(w)\}}{\displaystyle\sum_{J\subseteq C} f_J[h] \overline{w}_{(q\uparrow J)} u^{(q)}_J(w)},$$

where in the second line we use the *average cross-product ratio* $R^{(q)}_{J\bullet}(w) = k_J^{-1} \sum_{\substack{K\subseteq(C\setminus B)\setminus J \\ |K|=k}} R^{(q)}_{JK}(w)$, averaged over all subsets $K\subseteq C\setminus B$ of size $k$ not intersecting $J$), and where $u^{(q)}_J(w) = \sum_{(b-q)\subseteq C\setminus J} w_{(b-q)}$. In the third line we change the index of summation from $[j]$ to all subsets $J\subseteq C$ of size $|J|=p$, writing $f_J[h]$ for the number of permutations of the elements of $J$ for which $[j]\le[h]$, which number may equal zero for $J$ with all $[j]\not\le[h]$. Note that $\sum_{\substack{J\subseteq C \\ |J|=p}} f_J[h] = p!\binom{[h]}{p}$. The average value of the $w_{(q\uparrow j)}$ over subscripts with $f_J[h]>0$ we write as $\overline{w}_{(q\uparrow J)} = \overline{w}_{(q\uparrow J)}[h] = \sum_{[j]:\{j\}=J} I\{[j]\le[h]\} w_{(q\uparrow j)} / f_J[h]$, setting it equal to zero if $f_J[h]=0$. Thus the cusum in (4.2) is positive if and only if

$$\frac{\sum_{J\subseteq C} f_J[h]\overline{w}_{(q\uparrow J)}u_J^{(q)}(w)\{k_J R_{J\bullet}^{(q)}(w)\}}{\sum_{J\subseteq C} f_J[h]\overline{w}_{(q\uparrow J)}u_J^{(q)}(w)} > \frac{\binom{c-b-p+q}{k}\binom{c-p}{b-q}\binom{c-k}{b}}{\binom{c-k-p}{b-q}\binom{c}{b}}.$$ After manipulating factorials, the right-hand quantity simplifies to the two equivalent expressions $\dfrac{\binom{c-k}{b}\binom{c-p}{k}}{\binom{c}{b}} = \dfrac{\binom{c-b}{k}\binom{c-k}{p}}{\binom{c}{p}}$.

Thus the cusum in (4.2) is positive if and only if

$$\frac{\sum_{J\subseteq C} f_J[h]\overline{w}_{(q\uparrow J)}u_J^{(q)}(w)\{k_J R_{J\bullet}^{(q)}(w)\}}{\sum_{J\subseteq C} f_J[h]\overline{w}_{(q\uparrow J)}u_J^{(q)}(w)} > \frac{\binom{c-b}{k}\binom{c-k}{p}}{\binom{c}{p}}. \tag{4.4}$$

In words, cusum positivity is equivalent to a certain weighted average (over all subsets $J$ of size $p$) of $k_J = \binom{c-b-p+q'(J)}{k}$ times the average cross-product ratio being greater than the right-hand side of (4.4), which is independent of $q[i]$ and only depends on the constants $c$, $b$, $k$, and $p$.

It is instructive now to consider the special case in which all the odds parameters in $w$ are equal and without loss of generality we may take them all equal to one. In this case, all the $w$-product averages $\overline{w}_{(q\uparrow J)}$ and the average cross-product ratios $R_{J\bullet}^{(q)}(w)$ are equal one and the $u_J^{(q)}(w)$ are all equal, so the left-hand side of (4.4) is just the average value of $k_J$ over all $[j]\leq [h]$. Here we have the following result.

**Lemma 4.1.** *For any $p=1,\ldots,c-k$, any multi-index cusum superscript, and any contributing cusum subscript $h=(h_1,\ldots,h_p)$ with $1\leq h_\alpha \leq c-1$ for $\alpha = 1,\ldots p$, we have*

$$\frac{\sum_{\substack{J\subseteq C \\ |J|=p}} f_J[h]k_J}{p!\binom{[h]}{p}} = \underset{[j]\leq[h]}{\operatorname{avg}}\binom{c-b-p+q'[j]}{k} > \frac{\binom{c-k}{b}\binom{c-p}{k}}{\binom{c}{b}} = \frac{\binom{c-b}{k}\binom{c-k}{p}}{\binom{c}{p}}. \tag{4.5}$$

Thus the cusums in (4.4) are all positive when $w_1 = \cdots = w_c$. In words, the number of $K$ subsets not intersecting an $F$-component subscript $[j]$, when averaged over all $[j]$ involved in any contributing cusum subscript $[h]$, is bounded from below by the right-hand side of (4.5), independent of the cusum superscript $[i]$. A complete proof of Lemma 4.1 appears in the Appendix. Because the inequality in (4.5) is strict, cusum positivity holds for other odds parameters $w_1 \geq \cdots \geq w_c$ not too dissimilar one from another. We might say, then, that the behavior of the number of $K$-subsets involved in the cumulative summations of $F$-components and, in particular, the behavior of $q'[j]$ over $[j] \leq [h]$ described in (4.5) provides heuristic reason #1 for cusum positivity.

Next, consider the special case $p=1$, but for general $w$, wherein the single-index $F$-components have the following property.

**Lemma 4.2.** *In the single-index case $p=1$, for any superscript $[i]$, we have $F_{[j]}^{[i]}(a^+) > 0$ for $j=1,\ldots,b$ and $F_{[j]}^{[i]}(a^+) < 0$ for $j=b+1,\ldots,c$.*

Levin and Leu (2013) used such a *monotone sign property* and the fact that $S_{[c]}^{[i]}(a^+) = F(a^+) = 0$ to deduce the positive cusum property of Lemma 4.1, i.e., $S_{[h]}^{[i]}(a^+) > 0$ for the single-index cusums $h=1,\ldots,c-1$. For suppose it were the case that $S_{[h]}^{[i]}(a^+) \leq 0$ for some $h < c$. Then one (or more) of the preceding single-index $F$-components must be strictly negative; but then, by Lemma 4.2, all subsequent single-index $F$-components must also be strictly negative, in which case $F(a^+)$ would be negative, contradiction.

Levin and Leu (2013a) provided a partial proof of Lemma 4.2, namely, when $j \leq b$ with $q[i]=1$ (i.e., $i \leq b$) and when $j > b$ with $q[i]=0$ (i.e., $i > b$). In those two cases, the number of subsets $K$ in the bracketed sum on the left-hand side of (4.1) exactly matches the binominal coefficient on the right-hand side, such that for each $j$ and each $K$ in $(C \setminus B) \setminus \{j\}$ with $|K|=k$, a rigorous counting argument could be given for the respective positivity or negativity of the pairwise differences

$$\frac{\sum_{(b-q) \subseteq (C \setminus K) \setminus \{j\}} w_{(b-q)}}{\sum_{(b) \subseteq C \setminus K} w_{(b)}} - \frac{\sum_{(b-q) \subseteq C \setminus \{j\}} w_{(b-q)}}{\sum_{(b) \subseteq C} w_{(b)}}.$$ That counting argument, however, could not be extended to the

other two cases, namely, when $j \leq b$ with $q[i]=0$ (i.e., $i > b$) and when $j > b$ with $q[i]=1$ (i.e., $i \leq b$), because in those cases, the sums on the left-hand side of (4.1) are *excessive*, i.e., they have more terms than the binomial coefficient on the right-hand side of (4.1), or are *deficient*, i.e., have fewer terms, respectively, but the respective pairwise inequalities go in the wrong directions. So the pair-matching strategy fails, which is why Levin and Leu (2013) resorted to computer-assisted symbolic manipulation to confirm those cases, albeit only for $c \leq 7$.

Our purpose in revisiting the single-index case here is to view the monotone sign property of Lemma 4.2 in a new light that will provide another heuristic reason for the general PCP, and along the way, to give a new proof for the non-matched cases which is presented in Levin and Leu (2020). The development illustrates certain bounds on the average cross-product ratios in (4.4) leading to the desired result. Those properties are stated in the next lemma, which is equivalent to Lemma 4.2.

**Lemma 4.3.** *For $p=1$, $q[i]=0$, and $J=\{j\}$ for $j=1,\ldots,b$, the average cross-product ratios $R_{J\bullet}^{(q)}(w)$ are bounded from below by $1-(k/c)<1$. For $p=1$, $q[i]=1$, and $J=\{j\}$ for $j=b+1,\ldots,c$, the average cross-product ratios $R_{J\bullet}^{(q)}(w)$ are bounded from above by $\dfrac{1-(k/c)}{1-\{k/(c-b)\}}>1$.*

Parallel to how the PCP follows from Lemma 4.2, the PCP can be deduced from Lemma 4.3 as follows. Applying the same re-expressions leading from (4.2) to (4.4), when there are no restrictions to $F$-component subscripts $[j]\leq[h]$, we have the equality

$$\frac{\sum_{J\subseteq C}\overline{w}_{(q\uparrow J)}u_J^{(q)}(w)\{k_J R_{J\bullet}^{(q)}(w)\}}{\sum_{J\subseteq C}\overline{w}_{(q\uparrow J)}u_J^{(q)}(w)} = \frac{\binom{c-b}{k}\binom{c-k}{p}}{\binom{c}{p}},$$

corresponding to the constraint $S_{[c]}^{[i]}(a^+)=F(a^+)=0$. Furthermore, because the frequencies $f_J[h]$ can only equal one or zero when $p=1$, the equality is equivalent to

$$\frac{\sum_{J\subseteq C}f_J[h]\overline{w}_{(q\uparrow J)}u_J^{(q)}(w)}{\sum_{J\subseteq C}\overline{w}_{(q\uparrow J)}u_J^{(q)}(w)}Avg_+ + \frac{\sum_{J\subseteq C}(1-f_J[h])\overline{w}_{(q\uparrow J)}u_J^{(q)}(w)}{\sum_{J\subseteq C}\overline{w}_{(q\uparrow J)}u_J^{(q)}(w)}Avg_- = \frac{\binom{c-b}{k}\binom{c-k}{p}}{\binom{c}{p}}, \quad (4.6)$$

where

$$Avg_+ = \frac{\sum_{J\subseteq C}f_J[h]\overline{w}_{(q\uparrow J)}u_J^{(q)}(w)\{k_J R_{J\bullet}^{(q)}(w)\}}{\sum_{J\subseteq C}f_J[h]\overline{w}_{(q\uparrow J)}u_J^{(q)}(w)} \quad \text{and} \quad Avg_- = \frac{\sum_{J\subseteq C}(1-f_J[h])\overline{w}_{(q\uparrow J)}u_J^{(q)}(w)\{k_J R_{J\bullet}^{(q)}(w)\}}{\sum_{J\subseteq C}(1-f_J[h])\overline{w}_{(q\uparrow J)}u_J^{(q)}(w)}.$$

Therefore, the PCP property (4.4) can be established either by showing that $Avg_+$ is greater than the overall average in (4.6) or by showing that $Avg_-$ is less than the overall average. Now, when $j\leq b$, the number of subsets $K$ not intersecting $J$ equals $k_J=\binom{c-b-p+q'(J)}{k}=\binom{c-b}{k}$ with $p=1$, and because $\binom{c-k}{p}/\binom{c}{p}=1-(k/c)$, Lemma 4.3 implies that $k_J R_{J\bullet}^{(q)}(w)$ exceeds $\binom{c-b}{k}\binom{c-k}{p}/\binom{c}{p}$. Therefore, for any cusum subscript with $h\leq b$, the weighted average of those terms in $Avg_+$ also exceeds the overall average in (4.6). Similarly, when $j>b$, the number of subsets $K$ not intersecting $J$ equals $k_J=\binom{c-b-1}{k}$, and because $\binom{c-b}{k}\binom{c-k}{p}/\binom{c-b-1}{k}\binom{c}{p}=\dfrac{1-(k/c)}{1-\{k/(c-b)\}}$, Lemma 4.3 implies that $k_J R_{J\bullet}^{(q)}(w)$ is less than $\binom{c-b}{k}\binom{c-k}{p}/\binom{c}{p}$. Therefore, for any cusum subscript with $h>b$, the weighted average of those terms in $Avg_-$ is also less than the overall average in (4.6). This establishes the PCP for the single-index case.

In summary, we might say that the lower and upper bounds to the average cross-product ratios, working in conjunction with the partitioning induced by the terms included or excluded by summing $F$-components over $[j] \leq [h]$, provides heuristic reason #2 for cusum positivity.

In future work we hope to provide a general proof of Lemma 2.1 by combining the stochastic inequalities of Lemma 4.1 with bounds on the average cross-product ratios such as contained in Lemma 4.3. Here we close with an illustration of what we mean.

**Example 4.1.** Consider double-index cusums ($p=2$) with $q=q[i]=1$ and $[h]=(1,h_2)$ for $2 \leq h_2 \leq c-1$ and $b \geq 2$. Levin and Leu (2013) relied on the positivity of such cusums to help demonstrate the key inequality $F(a)>0$. To illustrate how inequality (4.4) dovetails with bounds on the average cross-product ratio to create multi-index cusum positivity, we use special odds parameters called *simplex boundary vectors* of the form $W^{(g)}(\omega) = (\omega,...,\omega,1,...,1)$ where, for $g=1,...,c-1$, the first $g$ components are equal to a given constant $\omega \geq 1$. For an explanation of this terminology and the rationale for considering simplex boundary vectors, see Levin and Leu (2020). Because the first $g$ components of these vectors are equal, the values of the average cross-product ratios $R_{J\bullet}^{(q)}(W^{(g)}(\omega))$ and other terms in equation (4.4) take on only a few distinct values, making it easier to highlight the relations. For present purposes, we only need to consider values of $g=b$, but we'll state results that apply for $g \leq b$. We utilize limiting values of $R_{J\bullet}^{(q)}(W^{(g)}(\omega))$ as $\omega$ grows large and write these as $R_{J\bullet}^{(q,g)}(\infty) = \lim_{\omega \to \infty} R_{J\bullet}^{(q)}(W^{(g)}(\omega))$. Such limiting values play an important role in establishing Lemma 4.3 in the single-index case too; see Levin and Leu (2020). Using similar methods, it is straightforward to derive result (4.7) below. The various properties of the average cross-product ratios and weighted averages stated below have been observed to hold empirically in all cases examined, and we surmise they hold in general.

We begin by noting that with $J = \{1, j_2\}$, $f_J[h]=1$ and $\overline{w}_{(q\uparrow J)} = \omega$ for $j_2 = 2,...,c-1$; because these values are constant, they have no effect on the weighted averages in (4.4) upon normalization of the weights. Also, $k_J = \binom{c-b}{k}$ for $j_2 = 2,...,b$ and $k_J = \binom{c-b-1}{k}$ for $j_2 = b+1,...,c-1$. Next, by considering the leading powers of $\omega$ in each of the four factors of the cross-product ratios, taking the limit as $\omega \to \infty$, and simplifying the binomial coefficients, we find

$$R_{J\bullet}^{(q,g)}(\infty) = \begin{cases} (1-k/c)\{1-k/(c-1)\} & \text{for } 2 \leq j_2 \leq g \\ \\ \dfrac{(1-k/c)\{1-k/(c-1)\}}{1-k/(c-g)} & \text{for } g+1 \leq j_2 \leq c-1. \end{cases} \quad (4.7)$$

As a function of $\omega$, $R_{J\bullet}^{(q)}(W^{(g)}(\omega))$ has what we call the *eventually-decreasing property*, namely, that regardless of whether $R_{J\bullet}^{(q)}(W^{(g)}(\omega))$ initially increases or decreases from 1 at $\omega=1$, it *eventually decreases* to $R_{J\bullet}^{(q,g)}(\infty)$ as $\omega \to \infty$. We note that for certain values of $b$, $c$, and $k$, when $j_2 > b$, the limiting value $R_{J\bullet}^{(q,g)}(\infty)$ in the second line of (4.7) may be greater than 1 and for such values it is thus not true that $R_{J\bullet}^{(q)}(W^{(g)}(\omega)) > R_{J\bullet}^{(q,g)}(\infty)$ *for all* $\omega$. However, it will suffice for present purposes for the inequality to hold *for sufficiently large* $\omega$.

Next, for an *arbitrary* odds parameter $w = (w_1, w_2, ..., 1)$ with $w_1 \geq \cdots \geq w_{c-1} \geq 1$, in the case under consideration, the weighted average of $k_J R_{J\bullet}^{(q)}(w)$ in (4.4) equals

$$Wtd\ Avg(w) = \frac{\sum_{j_2=2}^{h_2} u_J^{(q)}(w) k_J R_{J\bullet}^{(q)}(w)}{\sum_{j_2=2}^{h_2} u_J^{(q)}(w)}$$

with $u_J^{(q)}(w) = \sum_{(b-1) \subseteq C \setminus J} w_{(b-1)}$. We find that for any given $w$ there exists an $\omega^* = \omega^*(w_1) \geq w_1$ that satisfies $Wtd\ Avg(w) \geq Wtd\ Avg(W^{(b)}(\omega^*(w_1)))$. That is to say, there is always a simplex boundary vector $W^{(g)}(\omega^*) = (\omega^*, ..., \omega^*, 1, ..., 1)$ with $g=b$ components equal to $\omega^*$ which produces a smaller weighted average than the one produced by the given $w$. This dramatically simplifies the analysis. In practice, we find $\omega^*(w_1) = w_1$ suffices for widely dispersed components of $w$, and otherwise small multiples of $w_1$ suffice in all cases checked, but the relevant point for analysis is that for any given weighted average, one can always produce a smaller one using a boundary simplex vector with $b$ components equal to some finite multiple of $w_1$.

Now let $A(\omega) = Wtd\ Avg(W^{(b)}(\omega))$ for any $\omega \geq 1$. We observe that $A(\omega)$ also has the eventually-decreasing property: whether $A(\omega)$ initially increases or decreases from $A(1)$, for sufficiently large $\omega$ it eventually decreases to a finite limiting value as $\omega \to \infty$. In addition, we find that the derivative $(d/d\omega)A(\omega)$ *has at most two sign changes* on the interval $1 < \omega < \infty$. If the derivative has no sign changes, then $A(\omega)$ decreases monotonically from 1 to a limiting value as $\omega$ increases from 1 to $\infty$. If the derivative has one sign change, then $A(\omega)$ initially increases, reaches a maximum value, and eventually decreases to a limiting value as $\omega \to \infty$. If the derivative has two sign changes, then $A(\omega)$ initially decreases, rapidly reaching a local minimum value, then increases to a local maximum, and eventually decreases to a limiting value as $\omega \to \infty$. *In such cases, however, the initial local minimum always exceeds the limiting value, which is thus a global minimum for $A(\omega)$.* We think of this as a *quasi-unimodality property*, because if $(d/d\omega)A(\omega)$ has only one sign change, or none at all, then $A(\omega)$ is actually unimodal, which means it cannot initially decrease *below* the limiting value and then increase in order to eventually decrease again toward the limiting value.

When $(d/d\omega)A(\omega)$ has two sign changes, such "down-up-down" behavior will not invalidate the argument below because it will suffice for the initial local minimum to still exceed the limiting value. The example in Figure 1 below exhibits quasi-unimodality in the weighted average.

We now show that $A(\omega) \geq \binom{c-b}{k}\binom{c-k}{p}/\binom{c}{p}$ for *all* $\omega \geq 1$, and in particular $\omega = \omega^*(w_1)$, so that $Wtd\ Avg(w) \geq A(\omega^*(w_1))$ satisfies inequality (4.4). Thanks to the special structure of $W^{(b)}(\omega)$, the average cross-product ratios have two sets of common values, namely, $R^{(q)}_{\{1,2\}\bullet}(W^{(b)}(\omega)) = \cdots = R^{(q)}_{\{1,b\}\bullet}(W^{(b)}(\omega))$ and $R^{(q)}_{\{1,b+1\}\bullet}(W^{(b)}(\omega)) = \cdots = R^{(q)}_{\{1,c-1\}\bullet}(W^{(b)}(\omega))$, and similarly the weight components $u_j^{(q)}(w)$ have two sets of common values, $u^{(q)}_{\{1,2\}}(W^{(b)}(\omega)) = \cdots = u^{(q)}_{\{1,b\}}(W^{(b)}(\omega))$ and $u^{(q)}_{\{1,b+1\}}(W^{(b)}(\omega)) = \cdots = u^{(q)}_{\{1,c-1\}}(W^{(b)}(\omega))$. It follows that if $h_2 \leq b$, all the pertinent weight components are equal and $A(\omega)$ reduces simply to $A(\omega) = \binom{c-b}{k} R^{(q)}_{\{1,b\}\bullet}(W^{(b)}(\omega))$. Then by (4.7) and the quasi-unimodality property, we already have $A(\omega) > \binom{c-b}{k}(1-k/c)\{1-k/(c-1)\} = \binom{c-b}{k}\binom{c-k}{p}/\binom{c}{p}$.

When $b < h_2 \leq c-1$,

$$A(\omega) = U_b(\omega)\binom{c-b}{k}R^{(q)}_{\{1,b\}\bullet}(W^{(b)}(\omega)) + U_{b+1}(\omega)\binom{c-b-1}{k}R^{(q)}_{\{1,b+1\}\bullet}(W^{(b)}(\omega)),$$

where we have written the weights as $U_b(\omega) = \dfrac{(b-1)u^{(q)}_{\{1,b\}}(\omega)}{(b-1)u^{(q)}_{\{1,b\}}(\omega) + (h_2-b)u^{(q)}_{\{1,b+1\}}(\omega)} = 1 - U_{b+1}(\omega)$. Now, crucially, we have the initial value $A(1) = \dfrac{b-1}{h_2-1}\binom{c-b}{k} + \dfrac{h_2-b}{h_2-1}\binom{c-b-1}{k} > \binom{c-b}{k}\binom{c-k}{2}/\binom{c}{2}$ at $\omega = 1$, by Lemma 4.1. Also, for sufficiently large $\omega$,

$$A(\omega) > U_b(\omega)\binom{c-b}{k}R^{(q,b)}_{\{1,b\}\bullet}(\infty) + U_{b+1}(\omega)\binom{c-b-1}{k}R^{(q,b)}_{\{1,b+1\}\bullet}(\infty)$$

$$= U_b(\omega)\binom{c-b}{k}(1-k/c)\{1-k/(c-1)\} + U_{b+1}(\omega)\binom{c-b-1}{k}\dfrac{(1-k/c)\{1-k/(c-1)\}}{1-k/(c-b)}.$$

But $\binom{c-b-1}{k}\dfrac{(1-k/c)\{1-k/(c-1)\}}{1-k/(c-b)} = \binom{c-b}{k}\binom{c-k}{2}/\binom{c}{2}$, so $A(\omega) > \binom{c-b}{k}\binom{c-k}{2}/\binom{c}{2}$ for sufficiently large $\omega$. We conclude from the quasi-unimodality property that the inequality holds *for all* $\omega \geq 1$, as required.

Thus (4.4) holds and we have the positive cusum property for double-cusums of the form $S^{(b-1,b+1)}_{(1,h_2)}(a^+)$, for example, with $h_2 = 2,\ldots,c-1$, as required in Levin and Leu (2013).

# Figure 1

Graphs of $\binom{c-b}{k} R^{(q)}_{\{1,b\}\bullet}(W^{(b)}(\omega))$ (top), $\binom{c-b-1}{k} R^{(q)}_{\{1,b+1\}\bullet}(W^{(b)}(\omega))$ (bottom), and *Wtd Avg* (middle)

where $b=g=25$, $c=48$, and $k=22$

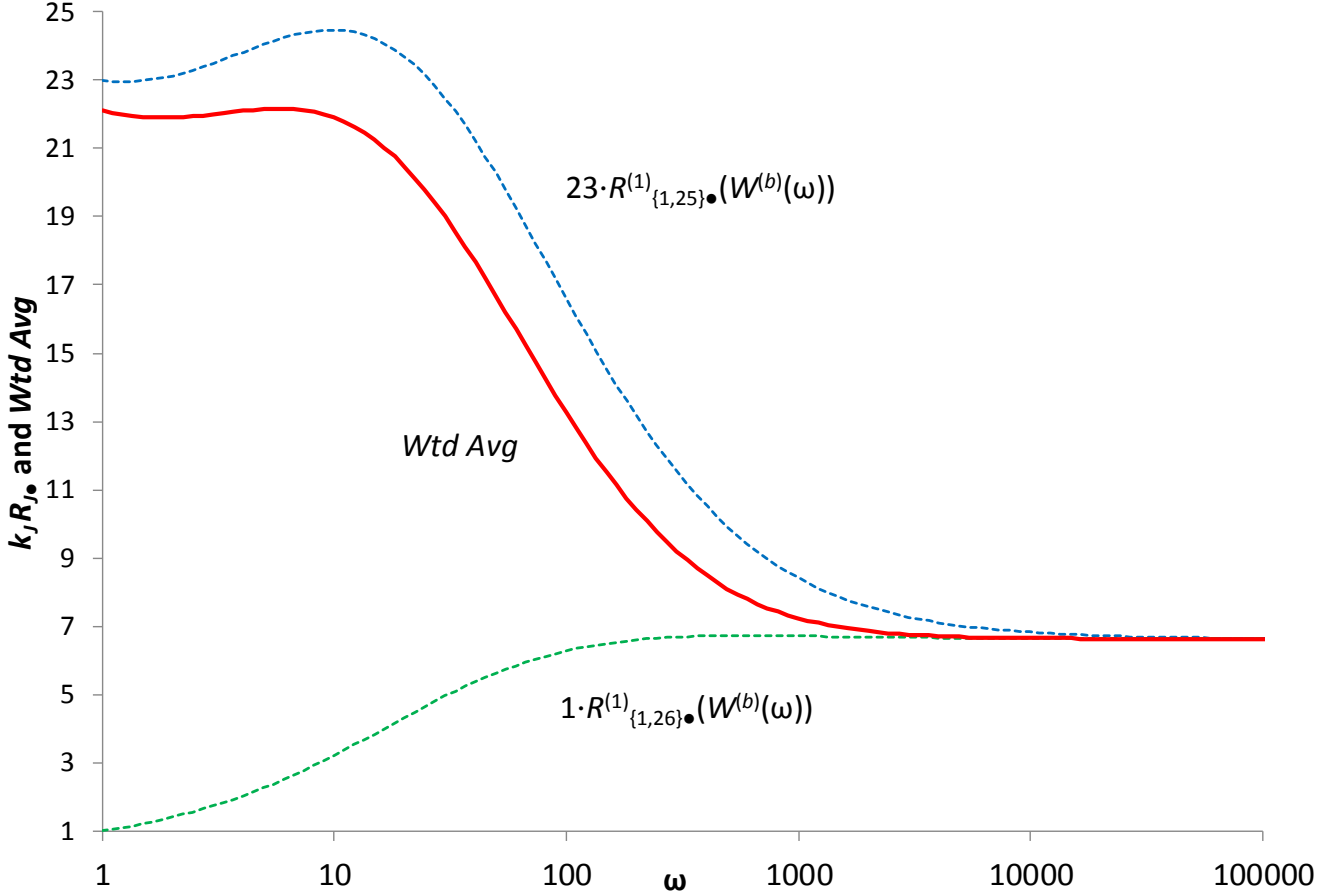

**ACKNOWLEDGMENT**

We thank the editors and referees for their generous efforts on our behalf.

**APPENDIX: Proof of Lemma 4.1**

We prove Lemma 4.1 by identifying all $[j] \leq [h]$ as the realizable values of certain random variables. To avoid confusion of notation, in this appendix we will only denote the subset of elements in subscript $[j]$ as $\{j\}$, rather than by $J$ as was done in the main text, which will then allow us instead to use $J = (J_1, ..., J_p)$ to denote the following random variables. Let $J_1$ be uniformly distributed over the set $\{1, ..., h_1\}$; then, given $J_1 = j_1$, let $J_2$ be uniformly distributed over the set $\{1, ..., h_1\} \setminus \{j_1\}$; and so on, until given $J_1 = j_1, ..., J_{p-1} = j_{p-1}$, let $J_p$ be uniformly distributed over the set $\{1,...,h_p\} \setminus \{j_1,...,j_{p-1}\}$. Also, define the indicator random variables $I_\alpha = I[J_\alpha \leq b]$ and their cumulative sums, $q'_\alpha = q'_\alpha[J] = \sum_{\ell=1}^{\alpha} I_\ell$ for $\alpha = 1,...,p$. Then the previously defined values of $q'[j]$ exhaust the possible realizations of the random variable $q'_p[J]$ and *vice versa*. Thus $\mathop{avg}\limits_{[j] \leq [h]} q'[j] = Eq'_p[J] = \sum_{\alpha=1}^{p} EI_\alpha = \sum_{\alpha=1}^{p} P[J_\alpha \leq b]$. Because the previously defined $q=q[i]$ will play no further role here, we henceforth drop the prime notation from $q'_p[J]$, writing $q_p[J]$ instead.

Note that in general, the components of $J$ are neither independent nor identically distributed. In fact, we claim that for each $\alpha = 1,...,p$, the random variable $q_\alpha[J]$ under the above probability model with the given cusum subscript $[h]$ is *stochastically strictly greater* (in the usual stochastic order) than under the probability model with equal subscript elements $[h'] = (h',...,h')$ for any $h' > h_\alpha$. That is, we claim that for any given $q$, $P_{[h]}[q_\alpha[J] > q] \geq P_{[h']}[q_\alpha[J] > q]$ with strict inequality for at least one $q$. The proof is by induction. The assertion is clear for $\alpha = 1$, because $P[J_1 \leq b] = b/h_1$ so that $q_1[J] = 1$ with probability $b/h_1$ else it equals zero. Thus there is the single support point $q=0$ with $P_{[h]}[q_1[J] > q] = b/h_1 > b/h'$ for any $h' > h_1$. Assume, then, that $P_{[h]}[q_\alpha[J] > q] \geq P_{[h']}[q_\alpha[J] > q]$ holds for each index up to a given $\alpha$ and consider $\alpha+1$. Then

$$P_{[h]}[q_{\alpha+1}[J] > q] = E_{[h]} P[q_{\alpha+1}[J] > q \mid q_\alpha[J]] = E_{[h]} P[I_{\alpha+1} > q - q_\alpha[J] \mid q_\alpha[J]]$$

$$= E_{[h]} \begin{cases} 1 & \text{if } q_\alpha[J] > q \\ \dfrac{b-q}{h_{\alpha+1} - \alpha} & \text{if } q_\alpha[J] = q \\ 0 & \text{if } q_\alpha[J] < q \end{cases}.$$

The bracketed function inside the expectation, call it $f_\alpha(\cdot \mid q, h_{\alpha+1})$, is non-decreasing as a function of the argument $q_\alpha[J]$, and strictly increases as $q_\alpha[J]$ varies from below, to equal, to greater than $q$ with positive probability. Under the inductive hypothesis, $q_\alpha[J]$ is stochastically greater under the probability model with $[h]$ than under the model with $[h']$, from which it follows that the expectation of $f_\alpha(q_\alpha[J] \mid q, h_{\alpha+1})$ under $[h]$ (which expectation depends only on $h_1,...,h_\alpha$) exceeds the corresponding expectation under $[h']$, i.e., $E_{[h]} f_\alpha(q_\alpha[J] \mid q, h_{\alpha+1}) > E_{[h']} f_\alpha(q_\alpha[J] \mid q, h_{\alpha+1})$. Furthermore, at the support point $q$ for the random variable $f_\alpha(q_\alpha[J] \mid q, h_{\alpha+1})$, the function $f_\alpha(\cdot \mid q, h_{\alpha+1}) > f_\alpha(\cdot \mid q, h')$ for any $h' > h_{\alpha+1}$; hence $E_{[h']} f_\alpha(q_\alpha[J] \mid q, h_{\alpha+1}) > E_{[h']} f_\alpha(q_\alpha[J] \mid q, h')$. Thus we have shown that $P_{[h]}[q_{\alpha+1}[J] > q] > P_{[h']}[q_{\alpha+1}[J] > q]$ for any $h' > h_\alpha$, which completes the induction. Intuitively, what the stochastic ordering means is that constraining the $J_\alpha$ to not exceed $h_\alpha$ increases the probability of smaller values of $J_\alpha$, leading to higher probabilities of larger values of $q_p[J]$, compared to when there are no constraints on the $J_\alpha$.

To complete the proof of (4.5) in Lemma 4.1, we note that $\binom{c-b-p+q_p[J]}{k}$ is a strictly increasing function of $q_p[J]$, so again by stochastic ordering, for any $[h'] = (h',\ldots,h')$ with $h' > h_p$,

$$\underset{[j] \leq [h]}{avg}\binom{c-b-p+q[j]}{k} = E_{[h]}\binom{c-b-p+q_p[J]}{k} > E_{[h']}\binom{c-b-p+q_p[J]}{k}.$$

In fact, let us consider $h'=c$. Even though $c$ is not a permissible value for a cusum subscript (because the cusum would then correspond to a smaller $p$), $h'=c$ is a permissible value in the probability model for generating the random variable $J$. Because there would then be no constraints whatsoever on the possible values of $J_\alpha$ from $1,\ldots,c$, the observed values of $J_1,\ldots,J_\alpha$ would be in one-to-one correspondence with a simple random sample without replacement of $p$ items from $C=\{1,\ldots,c\}$ and $q_p[J]$ would then have a standard hypergeometric distribution with parameters $b$, $p$, and $c$. We next demonstrate that the expected value of $\binom{c-b-p+q_p[J]}{k}$ under this simple random sampling model equals the value on the right-hand side of (4.5).

Let $M$ be an $H \times \binom{c-b}{k}$ matrix where $H = p!\binom{[h]}{p}$ with rows corresponding to the possible distinct values of $[j]$ and columns corresponding to the possible subsets $K \subseteq C \setminus B$ of size $|K|=k$, with entries in $M$ given by $I[K \cap \{j\} = \emptyset]$ indicating that subscript $[j]$ in a given row $\beta$ does not intersect subset $K$ in a given column $\gamma$. The binomial coefficients $\binom{c-b-p+q_p[j]}{k}$, which count the number of $K$ that do not intersect $[j]$, are simply the row margins of $M$ and the column margins are the number of $[j]$ that do not intersect a given $K$. Therefore

$$\underset{[j] \leq [c]}{avg}\binom{c-b-p+q[j]}{k} = \sum_{\beta=1}^{H}\sum_{\gamma=1}^{\binom{c-b}{k}} M_{\beta\gamma}/H = \sum_{\gamma=1}^{\binom{c-b}{k}}\left(\sum_{\beta=1}^{H} M_{\beta\gamma}/H\right) = \sum_{\gamma=1}^{\binom{c-b}{k}} P[K_\gamma \cap \{J\} = \emptyset].$$

But under simple random sampling without replacement, the probability of $J$ not intersecting any given $K \subseteq C \setminus B$ of size $|K|=k$ is constant and equal to $\binom{c-k}{p}/\binom{c}{p}$. Therefore

$$\underset{[j] \leq [c]}{avg}\binom{c-b-p+q[j]}{k} = \binom{c-b}{k}P[K \cap \{J\} = \emptyset] = \binom{c-b}{k}\binom{c-k}{p}/\binom{c}{p}.$$

Therefore we have shown that

$$\underset{[j]\leq[h]}{avg}\binom{c-b-p+q[j]}{k} > \binom{c-b}{k}\binom{c-k}{p}\bigg/\binom{c}{p} = \binom{c-k}{b}\binom{c-p}{k}\bigg/\binom{c}{b}.$$

This proves Lemma 4.1.